\documentclass[12pt]{amsart}
\usepackage{amssymb, amsmath, amsthm}
\usepackage[all]{xy}

\newcommand{\mr}[1]{\mathrm{#1}}
\newcommand{\mf}[1]{\mathfrak{#1}}
\newcommand{\mc}[1]{\mathcal{#1}}

\newcommand{\Z}{{\bf Z}}
\newcommand{\Q}{{\bf Q}}
\newcommand{\zp}{{\bf Z}_p}
\newcommand{\qp}{{\bf Q}_p}

\newcommand{\pr}[1]{Q_{#1}}
\newcommand{\numpr}[1]{|\pr{#1}|}

\DeclareMathOperator{\Gal}{Gal} 
 
\DeclareMathOperator{\Hom}{Hom} 

\def\Cl{\mathrm{Cl}}

\newcommand{\UT}[1]{\mathcal{O}_{#1}^{\times}}

\newtheorem{theorem}{Theorem}[section]

\newtheorem{lemma}[theorem]{Lemma}
\newtheorem{corollary}[theorem]{Corollary}

\theoremstyle{definition}

\newtheorem{example}[theorem]{Example}
\newtheorem{question}[theorem]{Question}

\theoremstyle{remark}
\newtheorem*{remark}{Remark}
\newtheorem*{remarks}{Remarks}
\newtheorem*{acknowledgments}{Acknowledgments}

\DeclareMathOperator{\rk}{rank} \DeclareMathOperator{\rank}{rank}
\DeclareMathOperator{\cork}{corank}

\DeclareMathOperator{\coker}{coker}

\renewcommand{\baselinestretch}{1.3}

\begin{document}

\renewcommand{\baselinestretch}{1}
\begin{abstract}
Consider the family of CM-fields which are pro-$p$ $p$-adic Lie
extensions of number fields of dimension at least two, which
contain the cyclotomic $\zp$-extension, and which are ramified at
only finitely many primes.  We show that the Galois groups of the
maximal unramified abelian pro-$p$ extensions of these fields are
not always pseudo-null as Iwasawa modules for the Iwasawa algebras
of the given $p$-adic Lie groups. The proof uses Kida's formula
for the growth of $\lambda$-invariants in cyclotomic
$\zp$-extensions of CM-fields. In fact, we give a new proof of
Kida's formula which includes a slight weakening of the usual $\mu
= 0$ assumption. This proof uses certain exact sequences involving
Iwasawa modules in procyclic extensions. These sequences are
derived in an appendix by the second author.
\end{abstract}

\renewcommand{\baselinestretch}{1.3}
\title[On the failure of pseudo-nullity]{On the failure of pseudo-nullity\\
of Iwasawa modules}
\author{Yoshitaka Hachimori and Romyar T. Sharifi}
\maketitle

\section{Introduction}
Let $p$ be a prime number.  Given a Galois extension $L$ of a
number field $F$, one may consider the inverse limit $X_L$ under
norm maps of the $p$-parts of class groups in intermediate number
fields in $L$. By class field theory, this is none other than the
Galois group of the maximal unramified abelian pro-$p$ extension
of $L$.  Setting $\mc{G} = \Gal(L/F)$, we let $\Lambda(\mc{G})$
denote the Iwasawa algebra $\zp[[\mc{G}]]$ of $\mc{G}$. By a
module for $\Lambda(\mc{G})$, we will always mean a compact left
$\Lambda(\mc{G})$-module.  In particular, $X_L$ may be given the
structure of a $\Lambda(\mc{G})$-module via conjugation, and we
shall be concerned in this article with its resulting structure.

Let $K$ denote the cyclotomic $\zp$-extension of $F$. We say that
a $p$-adic Lie extension $L/F$ is {\em admissible} if $\mc{G} =
\Gal(L/F)$ has dimension at least $2$, $L$ contains $K$, and $L/F$
is unramified outside a finite set of primes of $F$.  An
admissible extension will be said to be {\em strongly admissible}
if $\mc{G}$ is pro-$p$ and contains no elements of order $p$.  We
remark that any compact $p$-adic Lie group contains a pro-$p$
$p$-adic Lie group with no elements of order $p$ as an open
subgroup (cf.\ \cite[Corollary 8.34]{DdMS}).

R. Greenberg considered the situation in which $L$ is the
compositum of all $\zp$-extensions of $F$.  In this case,
$\Lambda(\mc{G})$ is isomorphic to a power series ring over $\zp$
in finitely many commuting variables.  Greenberg conjectured that
the annihilator of $X_L$ has height at least $2$ as a
$\Lambda(\mc{G})$-module, which is to say that $X_L$ is
$\Lambda(\mc{G})$-pseudo-null (see
\cite[Conjecture~3.5]{greenberg}).

In \cite{venj-str}, O. Venjakob generalized the notion of a
finitely generated pseudo-null $\Lambda(\mc{G})$-module to pro-$p$
groups $\mc{G}$ with no elements of order $p$.  In \cite{css},
this definition was further generalized to finitely generated
modules over arbitrary rings.  For $\mc{G}$ a compact $p$-adic Lie
group, a finitely generated $\Lambda(\mc{G})$-module $M$ is
pseudo-null if $\mr{Ext}^i_{\Lambda(\mc{G})}(M,\Lambda(\mc{G})) =
0$ for $i = 0,1$.
The following question became of interest.

\begin{question} \label{false}
    Let $L/F$ be an admissible $p$-adic Lie
    extension, and set $\mc{G} = \Gal(L/F)$.
    Is $X_L$ necessarily pseudo-null as a $\Lambda(\mc{G})$-module?
\end{question}

Set $\Gamma = \Gal(K/F)$.  In general, $X_K$ is known to be a
finitely generated torsion $\Lambda = \Lambda(\Gamma)$-module and
was conjectured by K.\ Iwasawa to have finite $\zp$-rank.  In
other words, Iwasawa conjectured that $X_K$ has trivial
$\mu$-invariant $\mu(X_K)$.  When $F/\Q$ is abelian, $\mu(X_K) =
0$ by a result of B.\ Ferrero and L.\ Washington \cite{FeWa}. In
the case that $F = \Q(\mu_p)$, the $\lambda$-invariant
$\lambda(X_K)$ is positive if and only if $p$ is an irregular
prime.  Thus, $X_K$ is often not pseudo-null as a
$\Lambda$-module, since pseudo-null $\Lambda$-modules are finite,
hence the reason for the dimension at least $2$ condition.

\begin{remark}
  The condition that $L$ contain $K$ is also necessary, as
  Greenberg (following Iwasawa) had constructed examples of
  extensions $L/F$ for which $X_L$ is not pseudo-null
  and $\mc{G} = \Gal(L/F)$ is
  free of arbitrarily large finite rank over $\zp$ but $L$ does
  not contain $K$ (unpublished).  In particular, these $X_L$ have
  nonzero $\mu$-invariants as $\Lambda(\mc{G})$-modules
  (cf.\ \cite[Section 3]{venj-str} for the general definition).
\end{remark}

In general, if $L/F$ is a admissible $p$-adic Lie extension, then
$X_L$ is a finitely generated torsion module over
$\Lambda(\mc{G})$ 
(cf.\ \cite[Theorem 7.14]{Ho} and Lemma \ref{torsion}). 
If, in addition, it
is finitely generated over $\Lambda(G)$, for $G = \Gal(L/K)$, then
the property of $X_L$ being $\Lambda(\mc{G})$-pseudo-null is
equivalent to its being $\Lambda(G)$-torsion (in an appropriate
sense; cf.\ Lemma \ref{ven-p-null}). We note that if Iwasawa's
$\mu$-invariant conjecture holds, then $X_L$ is indeed finitely
generated over $\Lambda(G)$ (again, see Lemma \ref{torsion}). So,
Question \ref{false} could very well be rephrased to ask if $X_L$
is a finitely generated $\Lambda(G)$-torsion module.

As we shall demonstrate in this paper, the answer to Question
\ref{false} is ``no,'' with counterexamples occurring frequently
for CM-fields $L$.  We remark that, until this point, no such
counterexamples had been known (or expected).  Note that if $L$ is
a CM-field, then complex conjugation provides a canonical
involution on $X_L$, and if $p$ is odd we obtain a canonical
decomposition $X_L = X_L^+ \oplus X_L^-$ into its plus and minus
one eigenspaces, respectively. In fact, we can compute the {\em
$\Lambda(G)$-rank} (see Section \ref{modules}) of $X_L^-$ for any
strongly admissible $p$-adic Lie extension $L/F$ of CM-fields with
$\mu(X_K^-) = 0$. This rank is zero if and only if $X_L^-$ is
$\Lambda(G)$-torsion, or equivalently,
$\Lambda(\mc{G})$-pseudo-null.

\begin{theorem} \label{main}
    Let $p$ be an odd prime and
    $L/F$ be a strongly admissible $p$-adic Lie extension of CM-fields,
    with Galois group $\mc{G}$.
    Assume that $\mu(X_K^-) = 0$.
    Let $\pr{L/K}$ be the set of primes in the maximal real subfield
    $K^+$ of $K$ that split in $K$, ramify in $L^+$, and do not divide $p$.
    Let $\delta$ be $1$ or $0$ depending upon whether $F$ contains the
    $p$th roots of unity or not, respectively.
    Then we have
    $$
    \rk_{\Lambda(G)}(X_L^-)= \lambda(X_K^-)-\delta+\numpr{L/K}.
    $$
    In particular, $X_L^-$ is not pseudo-null over
    $\Lambda(\mc{G})$ if
    and only if
    $\lambda(X_K^-) - \delta + \numpr{L/K} \geq 1$.
\end{theorem}

The proof of Theorem \ref{main} uses a formula of Y.\ Kida's
\cite{Ki} for $\lambda$-invariants in CM-extensions of cyclotomic
$\zp$-extensions of number fields.  We will also give another
proof of Kida's formula.

Note that Theorem \ref{main} provides counterexamples to a
positive answer to Question \ref{false} even in the case that $F =
\Q(\mu_p)$ and $L$ is a $\zp$-extension of $K$ with complex
multiplication which is unramified outside $p$ (see Example
\ref{basicex}). For instance, the smallest prime $p$ for which the
$\zp$-rank of $X_K^-$ is (at least) 2 is $p = 157$.
By Kummer duality, there are two $\zp$-extensions $L$ of $K$ that
are Galois over $\Q$ and for which $X_L^-$ is not pseudo-null as a
$\Lambda(\mc{G})$-module. Other examples occur for $p = 353, 379,
467, 491,$ and so on.

On the other hand, some mild evidence for a positive answer to
Question \ref{false} is given in \cite{massey} and \cite{paireis}
in the case that $F = \Q(\mu_p)$ and $L$ is a $\zp$-extension of
$K$ which is unramified outside $p$ and defined via Kummer theory
by a sequence of cyclotomic $p$-units in $\Q(\mu_{p^{\infty}})^+$.
In particular, such $L$ are not CM-fields.  For instance, it is
shown in \cite{paireis} that, for the $\zp$-extension $L =
K(p^{1/p^{\infty}})$ of $K$ with $F = \Q(\mu_p)$, the
$\Lambda(\mc{G})$-module $X_L$ is pseudo-null for all $p < 1000$.
So, even with counterexamples to pseudo-nullity, there remains the
question of finding a natural class of extensions $L/F$ over which
$X_L$ is pseudo-null as a $\Lambda(\mc{G})$-module.

We describe a couple of possibilities for such a class, though
this description is tangential to the rest of the paper. First,
consider an algebraic variety $Z$ over $F$, and form, for some $i
\ge 0$ and $r \in \Z$, the cohomology group
$H^i_{\text{\'et}}(Z,\qp(r))$. In the spirit of Fontaine-Mazur
\cite{fm}, we say that $L/F$ {\em comes from algebraic geometry}
if $L$ lies in the fixed field of the representation of
$\Gal(\bar{\Q}/F)$ on such a cohomology group.  We mention the
following refinement of Question \ref{false}.

\begin{question} \label{possible}
    If $L/F$ is an admissible $p$-adic Lie extension which
    comes from algebraic geometry, then must $X_L$ be pseudo-null as
    a $\Lambda(\mc{G})$-module?
\end{question}

We feel that there is currently insufficient evidence for a
positive answer to this question to conjecture it in general.
However, it seems quite reasonable that it could hold, since we
restrict to a setting in which the size of $X_L$ might be
controlled by $p$-adic $L$-functions.  We believe that it is not
known if CM-fields arising as admissible $p$-adic Lie extensions
ever come from algebraic geometry, though it is generally expected
that they do not.

One might wish for a still larger class in Question
\ref{possible}, since even in the case that $F = \Q(\mu_p)$ and
$L$ is a $\zp$-extension of $K$ defined by a sequence of
cyclotomic $p$-units, the extension $L$ need not come from
algebraic geometry (if the Tate twist of $G$ is non-integral). So,
consider a tower of algebraic varieties $(Z_n)_{n \ge 0}$ defined
over $F$ such that the $Z_n$ are all Galois \'etale covers of $Z =
Z_0$ and the Galois group of the tower is a $p$-adic Lie group. We
then expand our class to contain those admissible $p$-adic Lie
extensions which lie in the fixed field of the Galois action on
$\displaystyle \lim_{\leftarrow}
H^i_{\text{\'et}}((Z_n)_{/\bar{\Q}},\qp(r))$ for some $i \ge 0$
and $r \in \Z$, in which the inverse limit is taken with respect
to trace maps (cf.\ \cite{ohta}).  Then any $L$ arising from
cyclotomic $p$-units can be recovered, for instance, from the
first cohomology groups with $\qp$-coefficients in the tower of
Fermat curves $x^{p^n} + y^{p^n} = z^{p^n}$ over $F = \Q(\mu_p)$
(cf.\ \cite[Corollary 1 of Theorem B]{iky}).

The organization and contents of this paper are as follows.   In
Section \ref{kida}, we give another proof of Kida's formula
(Theorem \ref{kida2}) which includes a slight weakening of the
usual assumption $\mu(X_K^-) = 0$.  We consider, in this formula,
any quotient $X_{L,T}$ of $X_L$ by the decomposition groups at
primes above $T$, for a finite set of primes $T$ of $F$. In
Section \ref{modules}, we discuss Iwasawa modules in $p$-adic Lie
extensions, elaborating on some of the definitions and remarks
given in this introduction as well as providing lemmas for later
use.  In Section \ref{general}, we prove the generalization of
Theorem \ref{main} to the case of $X_{L,T}$ (Theorem
\ref{second}).  In Section \ref{examples}, we provide, along with
a few remarks, specific examples of cases in which
$\Lambda(\mc{G})$-pseudo-nullity fails.  Finally, in Appendix
\ref{special}, the second author derives two exact sequences
involving the $G$-invariants and coinvariants of $X_L$ for quite
general procyclic extensions $L/K$ (Theorem \ref{comparison} and
Corollary \ref{zpext}).  These are used in the proof of Theorem
\ref{kida2}.

\begin{acknowledgments}
We wish to express our gratitude to John Coates for arranging for
our collaboration on this work, for his continued support
throughout,
and for his very helpful comments.  We also thank Susan Howson,
Manfred Kolster, Masato Kurihara, Kazuo Matsuno, Manabu Ozaki, and
Otmar Venjakob for their valuable advice.  
The first author was partially supported by Gakushuin University and 
the 21st Century COE program at the Graduate School of Mathematical Sciences of
the University of Tokyo.
The second author was supported by the Max Planck Institute for Mathematics.
\end{acknowledgments}

\section{Kida's formula} \label{kida}

In this section, we will give a proof of a mild generalization of
Kida's formula in which the condition on the $\mu$-invariant is
weakened.  For this, we use the exact sequences of Iwasawa modules
in cyclic extensions of Appendix \ref{special} (Theorem
\ref{comparison}). We let $p$ be an odd prime in this section.

Let $F$ denote a number field which is CM, and let $K$ denote its
cyclotomic $\Z_p$-extension. Consider a CM-field $L$ Galois over
$K$, and set $G = \Gal(L/K)$.  Let $T$ be a finite set of primes
of $F^+$
	(which we could just as well assume consists solely of primes above $p$
	in what follows)
, and for any algebraic extension $E/F^+$, let $T_E$ be
the set of primes above $T$. Define $P_K$ to be the set of primes
$v$ of $K$ lying above $p$. Fix a set $V_K^-$ consisting of one
prime of $K$ for each prime of $K^+$ that splits in $K$. Let
$T_K^- = T_K \cap V_K^-$ and $P_K^- = P_K \cap V_K^-$.
Furthermore, we let
\begin{equation} \label{primeset}
    Q_{L/K}^T =
    \{ v \in V_K^- - (T_K^- \cup P_K^-) \colon I_v \neq 1 \}
    \cup \{ v \in T_K^- \colon  G_v \neq 1 \},
\end{equation}
with $G_v$ and $I_v$ denoting the 
decomposition and inertia
groups
in $G$, respectively, at a chosen prime above $v$ in $L$.  That
is, $Q_{L/K}^T$ is in one-to-one correspondence with the set of
primes $u$ of $K^+$ such that $u$ splits in $K$ 
	but not completely in $L$ and either $u \in T_{K^+}$ or $u$ does not lie 
	above $p$.
For $v \in Q_{L/K}^T$, let $g_v(L/K) = [G:G_v]$.

Denoting the maximal real subfield of any $L$ as above by $L^+$,
we have a direct sum decomposition of any
$\zp[\Gal(L/L^+)]$-module $M$ into $(\pm 1)$-eigenspaces $M^{\pm}$
under complex conjugation. We shall be particularly interested in
(the minus parts of) two such Iwasawa modules.  That is, we define
$X_{L,T}$ to be the maximal unramified abelian pro-$p$ extension
of $L$ in which all primes in $T_L$ split completely, and we let
$\mc{U}_{L,T}$ denote the inverse limit of the $p$-completions of
the $T$-unit groups of the finite subextensions of $F$ in $L$.

We set
$$
    \delta = \begin{cases} 1 & \mr{if\ } \mu_p \subset F
    \\ 0 & \mr{if\ } \mu_p \not\subset F.
    \end{cases}
$$
For a finitely generated $\Lambda$-module $M$ and any choice of
pseudo-isomorphism
$$
    M(p) \to \bigoplus_{i=1}^N \Lambda/p^{r_i}\Lambda,
$$
where $M(p)$ denotes its $p$-power torsion subgroup,
set
$$
    \theta(M) = \max\,\{ r_i \mid 1 \le i \le N \}.
$$

The following remark should be kept in mind in what follows.

\begin{remark}
  Let $L/K$ be a finite Galois $p$-extension such that $L/F$ is
  Galois.  Set $\mc{G} = \Gal(L/F)$.
  Any $\Lambda(\mc{G})$-module $M$
  becomes a $\Lambda = \Lambda(\Gamma)$-module through a
  choice of subgroup of $\mc{G}$ lifting $\Gamma$, and the isomorphism class
  of $M$ as a $\Lambda$-module is independent of this choice.
\end{remark}

We have the following generalization of Kida's formula for the
behavior of $\lambda(X_L^-)$ in finite Galois $p$-extensions $L/K$
(such that $L/F$ is Galois) \cite{Ki} (see also the related
results of L.\ Kuz'min \cite[Appendix 2]{kuzmin} and in
\cite[Corollary 11.4.13]{nsw}, though for the latter we remark
that the formulas are not quite correct as written).  In the above
results, it is assumed that $\mu(X_{K,T}^-) = 0$, which implies
that $\mu(X_{L,T}^-) = 0$. (In these results, it is also assumed
that either $T$ is empty or consists of the primes above $p$.)  We
make the weaker assumption that $\theta(X_{L,T}^-) \le 1$. Of
course, Iwasawa conjectured that $\mu(X_{K,T}) = 0$ always holds.

\begin{theorem} \label{kida2}
    Let $p$ be an odd prime and $L$ a finite
    $p$-extension of $K$ which is Galois over $F$.
    Assume that $\theta(X_{L,T}^-) \le 1$.  Then
    \begin{equation} \label{kidasformula}
        \lambda(X_{L,T}^-)=[L:K](\lambda(X_{K,T}^-)
          -\delta+|Q_{L/K}^T|)+\delta
        -\sum_{v\in Q_{L/K}^T} g_v(L/K)
    \end{equation}
    and
    $$
        \mu(X_{L,T}^-) = [L:K]\mu(X_{K,T}^-).
    $$
\end{theorem}

\begin{proof}
    For notational convenience, we leave out superscript and
    subscript $T$'s in this proof.
    We claim that it suffices to
    demonstrate
    this result in the case that $[L:K] = p$.
    To see this, let $E$ be an intermediate field in $L/K$, and, if necessary,
    replace $F$ by a finite extension $F'$ of $F$ in $K$ such that
    $E/F'$ is Galois.
    Since $X_L^- \to X_E^-$ is surjective (as $G = G^+$), the fact that
    $\theta(X_L^-) \le 1$ implies that $\theta(X_E^-) \le 1$ as
    well.  Now we use induction on the degree of $[L:K]$ and assume that
    $L/E$ has degree $p$.  Then \eqref{kidasformula} for $L/K$
    amounts to a check of the formula
    $$
        \sum_{v \in Q_{E/K}} ([L:K]-p g_v(E/K))+ \sum_{w \in
        Q_{L/E}} (p-1) =\sum_{v \in Q_{L/K}} ([L:K]-g_v(L/K)),
    $$
    which we leave to the reader.
    Let us use $\mu'$ to denote $\mu$-invariants with respect to
    $\Gal(K/F')$.  By induction, we have
	\begin{align*}
		\mu(X_L^-) &= [F':F]^{-1}\mu'(X_L^-) = [F':F]^{-1}p\mu'(X_E^-)\\
		&= [F':F]^{-1}[L:K]\mu'(X_K^-) = [L:K]\mu(X_K^-).
	\end{align*}
    This proves the claim.

    So, let $G$
    be cyclic of order $p$.  We examine the minus parts of
    the sequences $\Gamma_{L/K}$ and $\Psi_{L/K}$ of Theorem \ref{comparison}.
    Since $\mc{U}_L^- \cong \zp(1)^{\delta}$ (see, for example,
    \cite[Theorem 11.3.11(ii)]{nsw}), we have
    $\hat{H}^0(\mc{U}_L^-) \cong \mu_p^{\delta}$ and
    $\hat{H}^{-1}(\mc{U}_L^-) = 0$.  Let $S_{L/K}$
    be the set of elements of $P_K^-$ with $v \notin T_K^-$
    and $I_v \neq 1$.  Since $G$ is cyclic of order $p$
    and $G = G^+$, the sequence $\Psi_{L/K}$ reduces to
    \begin{multline*}
    \Psi_{L/K}^-\colon\ \ 0 \to (G^{\otimes 2})^{\oplus |S_{L/K}|}
    \to \hat{H}^{0}(X_L^-) \otimes G \to \mu_p^{\delta} \to \\
    G^{\oplus |Q_{L/K}|+|S_{L/K}|} \to \hat{H}^{-1}(X_L^-) \to 0.
    \end{multline*}
    This yields that
    \begin{equation} \label{herbrand}
        h(X_L^-) = p^{\delta-|Q_{L/K}|},
    \end{equation}
    where we use $h(X_L^-)$ to denote the Herbrand quotient with respect to $G$.

    As for the minus-part of $\Gamma_{L/K}$, need from it only
    the rather well-known consequence that the restriction map
    $(X_L^-)_G \to X_K^-$ is a pseudo-isomorphism.  We have an
    exact sequence
    $$
        0 \to X_L^-(p) \to X_L^- \to Z \to 0,
    $$
    with $Z$ a $\Lambda[G]$-module which is free of finite rank
    over $\zp$.  Clearly, we have
    \begin{eqnarray} \label{lambda}
        \lambda(Z) = \lambda(X_L^-) & \mr{and} &
        \lambda(Z_G) = \lambda(X_K^-).
    \end{eqnarray}

    The sequence $\Psi_{L/K}^-$ also implies that
    $\mu(\hat{H}^i(X_L^-)) = 0$ 
    for all $i$.  As we shall see in
    Lemma \ref{mptrivial}(a), this and the
    the fact that $\theta(X_L^-) \le 1$ imply that $h(X_L^-(p)) =
    1$.  Therefore, we have
    \begin{equation} \label{extra}
        h(Z) = h(X_L^-).
    \end{equation}
    Now, we know by
    representation theory (as pointed out in \cite{Iw})
    that
    $$
        Z \cong \zp[G]^a \oplus I_G^b \oplus \zp^c,
    $$
    where $I_G$ denotes the augmentation ideal of $\zp[G]$.
    It follows that $\lambda(Z_G) = a+c$ and $h(Z) = p^{c-b}$.
    Applying \eqref{herbrand}, \eqref{lambda}, and \eqref{extra},
    we conclude that
    $$
        \lambda(X_L^-) = p(a+b) + c-b =
        p(\lambda(X_K^-)-\delta+|Q_{L/K}|)+\delta-|Q_{L/K}|,
    $$
    verifying \eqref{kidasformula}.

    As for the $\mu$-invariant, the fact
    that $(X_L^-)_G$ and $X_K^-$ are pseudo-isomorphic implies
    that $\mu(X_K^-) = \mu((X_L^-)_G)$.
    The result then follows from Lemma \ref{mptrivial}(b) below.
\end{proof}

To finish the proof of Theorem \ref{kida2}, we need the following
results on Herbrand quotients and $\mu$-invariants. Let $G$ be a
cyclic group of order $p$.  We use $\hat{H}^i(M)$ to denote the
$i$th Tate cohomology group of $M$ with respect to $G$.

\begin{lemma} \label{fptrivial}
    Let $M$ be a ${\bf F}_p[G]$-module
    for which $\hat{H}^i(M)$ is finite for all $i$.  Then
    the Herbrand quotient $h(M)$ is trivial.
\end{lemma}

\begin{proof}
    Note that ${\bf F}_p[G]$ has a filtration
    \begin{equation} \label{filtration}
        {\bf F}_p[G] = I_G^0 \supset I_G \supset \ldots \supset
        I_G^{p-1} = (N_G) \supset I_G^p = 0,
    \end{equation}
    where $I_G$ denotes the augmentation ideal and $N_G$ denotes
    the norm element in ${\bf F}_p[G]$.  Let $M[I_G^k]$ denote the
    submodule of
    elements of $M$ killed by $I_G^k$.  We then have exact
    sequences
    \begin{multline*}
        0 \to M[I_G^{k}]/(I_G M \cap M[I_G^{k}])
        \to M[I_G^{k+1}]/(I_G M \cap M[I_G^{k+1}])
        \xrightarrow{\phi_k} \\
        M^G/(I_G^{k+1} M)^G \to M^G/(I_G^k M)^G \to 0
    \end{multline*}
    for $k \ge 0$, and $\phi_0$ is
    simply the identity on $M^G/(I_G M)^G$.  Since the kernel and
    cokernel of $\phi_{k+1}$ are the domain and range of
    $\phi_k$, we conclude that the orders of the domain and
    range of $\phi_{p-2}$ are equal, if finite.  Since this
    finiteness is assumed, we conclude that $h(M) = 1$.
\end{proof}

In the following lemmas, we take $M$ to be a finitely generated
torsion $\Lambda$-module with a commuting $G$-action (i.e., a
finitely generated $\Lambda[G]$-module which is
$\Lambda$-torsion).  We use $N_G$ to denote the norm element in
$\zp[G]$.

\begin{lemma} \label{basiclemma}
    The following conditions on $M$
    are equivalent:
    \begin{enumerate}
    \item[(i)] $\mu(\hat{H}^0(M)) = 0$
    \item[(ii)] $\mu(\hat{H}^i(M)) = 0$ for all $i \in \Z$
    \item[(iii)] $\mu(M_G) = \mu(N_G M)$.
    \end{enumerate}
\end{lemma}

\begin{proof}
    The equivalence of (i) and (ii) follows immediately from
    $$
        0 \to M^G \to M \xrightarrow{\sigma-1} M \to M_G \to 0,
    $$
    for some generator $\sigma$ of $G$,
    \begin{equation} \label{2ndseq}
        0 \to \hat{H}^{-1}(M) \to M_G \xrightarrow{N_G} M^G \to \hat{H}^0(M)
        \to 0,
    \end{equation}
    and additivity of $\mu$-invariants in exact sequences.  The
    equivalence with (iii) also follows from \eqref{2ndseq}, as it
    implies that $\mu(M_G) = \mu(N_G M)$ if and only if
	$\mu(\hat{H}^{-1}(M)) = 0$.
\end{proof}

\begin{lemma} \label{mptrivial}
    Assume that $\mu(\hat{H}^0(M)) = 0$ and $\theta(M) \le 1$.
    Then we have:
    \begin{enumerate}
    \item[(a)] $h(M(p)) = 1$
    \item[(b)] $\mu(M) = p\mu(M_G)$.
    \end{enumerate}
\end{lemma}

\begin{proof}
    We first claim that it suffices to replace $M$ by
    $M(p)/pM(p)$.  For this,
    note that we have an exact sequence of $\Lambda[G]$-modules
    $$
        0 \to M(p) \to M \to Z \to 0,
    $$
    where $Z$ is free of finite rank over $\zp$.  Since $\mu(Z) = 0$,
    the same holds for its cohomology groups, and hence
    $$
        \mu(\hat{H}^0(M(p))) = \mu(\hat{H}^0(M)) = 0,
    $$
    which is to say that $\hat{H}^0(M(p))$ is finite. By Lemma
    \ref{basiclemma}, the groups $\hat{H}^i(M(p))$ are all finite.
    Since $\theta(M) \le 1$, we have that $M(p) \to M(p)/pM(p)$
    is a pseudo-isomorphism, and hence the $\hat{H}^i(M(p)/pM(p))$
    are finite as well.  As Herbrand quotients of finite modules
    are trivial, the claim is proven.

    With the claim proven and $M$ now assumed to be $p$-torsion, part (a)
    follows from Lemma \ref{fptrivial},
    using again the equivalence of (i) and (ii) in Lemma \ref{basiclemma}.
    For part (b), let
    $I_G$ denote the augmentation ideal in $\zp[G]$.
    Then $I_G^k M/I_G^{k+1} M$ is isomorphic to a
    quotient of $I_G^{k-1} M/I_G^k M$ for $k \ge 1$.  As in \eqref{filtration},
    we have
    $$
        \mu(I_G^{p-1} M/I_G^p M) = \mu(I_G^{p-1} M) = \mu(N_G M).
    $$
    Since $\mu(M_G) = \mu(N_G M)$
    by Lemma \ref{basiclemma}, it follows
    that
    $\mu(I_G^k M/I_G^{k+1} M) = \mu(M_G)$
    for $0 \le k \le p-1$.  We conclude that
    $$
        \mu(M) = \sum_{k=0}^{p-1} \mu(I_G^k M/I_G^{k+1} M) =
        p\mu(M_G).
    $$
\end{proof}

\begin{remark}
    In general, $h(M)$ can be nontrivial when $M$ is a $p$-power torsion
    $\Lambda$-module for which $h(M)$ is defined.
    For example, the principal $(\Lambda/p^2\Lambda)[G]$-ideal
    $M = ((g-1)-p(\gamma-1))$, for generators $g$ of $G$ and
    $\gamma$ of $\Gamma$, has $h(M) = p$.
\end{remark}

\section{Iwasawa modules} \label{modules}

In this section, we assemble a few definitions and
easily proven results regarding modules for
the Iwasawa algebras of $p$-adic Lie extensions.

Let $G$ be a compact $p$-adic Lie group.  We define a finitely
generated $\Lambda(G)$-module $M$ to be {\em torsion} (resp., {\em
pseudo-null}) if
$$
  \mr{Ext}^i_{\Lambda(G)}(M,\Lambda(G)) = 0
$$
for $i = 0$ (resp., for $i = 0,1$).  These definitions coincide
with those of Venjakob \cite[Sections 2-3]{venj-str} and J.
Coates, P. Schneider, and R. Sujatha in \cite[Section 1]{css}.
When $\Lambda(G)$ is an integral domain (for instance, if $G$ has
no elements of finite order), this definition of
$\Lambda(G)$-torsion reduces to the usual one.

For any pro-$p$ $p$-adic Lie group $G$ with no $p$-torsion, we let
$\mc{Q}(G)$ denote the skew fraction-field of $\Lambda(G)$. For
any finitely generated $\Lambda(G)$-module $M$, we define the {\em
$\Lambda(G)$-rank} of $M$ as a $\Lambda(G)$-module, denoted
$\rk_{\Lambda(G)} M$, to be the dimension of
$\mc{Q}(G)\otimes_{\Lambda(G)} M$ as a left $\mc{Q}(G)$-vector
space (see, for example, \cite[Section 2]{CH}). A finitely
generated $\Lambda(G)$-module $M$ is $\Lambda(G)$-torsion if and
only if it has trivial $\Lambda(G)$-rank.

We say that a compact $p$-adic Lie group $G$ is {\em uniform} (or,
{\em uniformly powerful}) if its commutator subgroup $[G,G]$ is
contained in the group $G^p$ generated by $p$th powers and the
$p$th power map induces isomorphisms on the successive graded
quotients in its lower central $p$-series (cf.\ \cite[Definition
4.1]{DdMS}).  It is known that any compact $p$-adic Lie group
contains an open normal subgroup which is uniform (cf.\
\cite[Corollary 8.34]{DdMS}).

\begin{lemma} \label{ven-p-null}
   Let $\mc{G}$ be a compact $p$-adic Lie group, and assume that
   $G$ is a closed normal subgroup with $\mc{G}/G \cong \zp$.
   Then a $\Lambda(\mc{G})$-module
   which is finitely generated over $\Lambda(G)$
   is pseudo-null if and only if it is torsion as a $\Lambda(G)$-module.
\end{lemma}

\begin{proof}
   This is proven in
   \cite[Proposition 5.4]{venjakob} (using \cite[Example 2.3]{venjakob})
   if $G$ is uniform and we have
   an inclusion of subgroups $[\mc{G},G] \le G^p$.
   We claim that $\mc{G}$ has
   an open subgroup $\mc{G}_1$ with this property, taking $G_1 =
   G \cap \mc{G}_1$.
   To see this, first let
   $\mc{G}_0$ be any open normal uniform subgroup of $\mc{G}$,
   and set $G_0 = G \cap \mc{G}_0$.  Fix $\gamma \in \mc{G}_0$
   with image generating $\mc{G}_0/G_0$, and let $\Gamma$ be the
   closed subgroup of $\mc{G}_0$ that $\gamma$ generates.
   We have a canonical isomorphism $\mc{G}_0 \cong G_0 \rtimes
   \Gamma$ (for the given action of $\Gamma$ on $G_0$).
   Let $G_1$ be an open normal uniform subgroup of $G_0$.
   Choose $n \ge 0$ such that $[\Gamma^{p^n},G_1] \le G_1^p$,
   and let $W$ be the closed normal subgroup of
   $G_1$ generated by $[\Gamma^{p^n},G_1]$.
   Let $\mc{G}_1 \cong G_1 \rtimes \Gamma^{p^n}$, an open
   subgroup of $\mc{G}$.
   Then $\mc{G}_1$ yields the claim, as
   $$
    [\mc{G}_1,G_1] \le [G_1,G_1] \cdot W \le G_1^p.
   $$

   The result now follows from the definitions of pseudo-nullity
   and torsion modules given above, since
   for any $\Lambda(\mc{G})$-module $M$, one has
   \cite[Lemma 2.3]{jannsen}
   $$
     \mr{Ext}^i_{\Lambda(\mc{G})}(M,\Lambda(\mc{G})) \cong
     \mr{Ext}^i_{\Lambda(\mc{G}_1)}(M,\Lambda(\mc{G}_1))
   $$
   for any $i \ge 0$,
   along with the corresponding fact replacing $\mc{G}$ by $G$
   and $\mc{G}_1$ by $G_1$.
 \end{proof}

We now consider the behavior of a certain sort of filtration on a
compact $p$-adic Lie group with respect to taking subgroups.

\begin{lemma} \label{well-known}
  Let $G$ be a compact $p$-adic Lie group, let $G_1$ be an open normal
  uniform pro-$p$ subgroup of $G$,
  and let $H$ be a closed subgroup of $G$.  For $n \ge
  1$, let $G_{n+1}$ denote the open normal subgroup of $G$
  which is the topological closure
  of $G_n^p[G_n,G_1]$ in $G$.
  Then $H$ is a $p$-adic Lie group of dimension $e \le \dim
  G$, and there exists a rational number $C$ such that
  $[H:H\cap G_n]=Cp^{ne}$
  for all sufficiently large $n$.
\end{lemma}

\begin{proof}
    The first statement is well-known, and as for the second
    statement,
    we begin by noting that
    $G_{n}={G_1^{p^{n-1}}}$
    since $G_1$ is uniform (cf.\ \cite[Theorem 3.6]{DdMS}).
    So, there exists an integer $c$ such that
    $H\cap G_{n}$ is
    uniform for all $n\geq c$
    (cf.\ \cite[\S 4 Exercise 14 (i)]{DdMS}).
    Let $H_0:=H\cap G_c$.
    Then we can take some $c'$ such that
    $H_0 \cap G_n=(H_0\cap G_{c'})^{p^{n-c'}}$ for all $n\geq c'$
    (cf.\ \cite[\S 4 Exercise 14 (ii)]{DdMS}).
    Setting $H_1=H_0\cap G_{c'}$ and defining
    $H_{n+1}$ to be the topological closure of $H_{n}^p[H_n,H_1]$ in $G$,
    we have
    $H_n= H_1^{p^{n-1}}=H_0\cap G_{n+c'-1}$
    since $H_1$ is uniform.
    As there exists some constant $C'$ such that
    $[H:H_0 \cap G_{n+c'-1}]=C'p^{ne}$ for all sufficiently
    large $n$ and $[H \cap G_n : H_n]$ is
    eventually constant, we have the result.
\end{proof}

We shall also
require the following asymptotic formula for $\Lambda(G)$-ranks
(cf.\ \cite[Theorem 2.22]{Ho}, or \cite[Theorem 1.10]{Ha} for
``adequate" $G$).  We say that a sequence $(a_n)_{n \ge 1}$ of
integers is (or ``equals'') $O(q^n)$ for some nonnegative integer
$q$ if $0 \le a_n \le Cq^n$ for some constant $C$ for all
sufficiently large $n$.

\begin{lemma}[Howson] \label{susan}
    Let $G$ be a pro-$p$ $p$-adic Lie group containing no
    elements of order $p$, and
    let $M$ be a finitely generated $\Lambda(G)$-module.
    Choose a sequence $G_n$ as in Lemma \ref{well-known}.
    Then $\rk_{\Lambda(G)} M=r$ if and only if
    $$
        \rk_{\zp} M_{G_n}=r[G:G_n]+O(p^{n(d-1)}).
    $$
\end{lemma}

Finally, we move away from the purely module-theoretic setting to
prove the following basically well-known consequence of Nakayama's
Lemma that was used in the introduction.

\begin{lemma} \label{torsion}
  Let $L/F$ be an admissible $p$-adic Lie extension with Galois
  group $\mc{G}$, and set $G = \Gal(L/K)$.  Then
  $X_{L,T}$ is a finitely generated 
  $\Lambda(\mc{G})$-module.
  Furthermore, if $L/K$ is strongly admissible and $\mu(X_{K,T}) = 0$,
  then $X_{L,T}$ is finitely generated as a $\Lambda(G)$-module.
\end{lemma}

\begin{proof}
  In the first part, 
  we may assume that $G$ is pro-$p$
  by passing to an open subgroup of $\mc{G}$.
  Note that $X_{K,T}$ is a finitely generated $\Lambda$-module.
  Furthermore, the kernel of $(X_{L,T})_G \to X_{K,T}$ is a finitely
  generated $\zp$-module,
  since this kernel 
	maps to the direct sum of the decomposition groups in $G$ at the ramified primes in $L/K$ and the primes in $T_K$,
	with kernel a quotient of $H_2(G,\zp)$.  
  Therefore, $(X_{L,T})_G$ is a finitely generated $\Lambda$-module,
  and it follows from Nakayama's Lemma (as in \cite{bh})
  that $X_{L,T}$ is a finitely
  generated $\Lambda(\mc{G})$-module.  If we know that
  $\mu(X_{K,T}) = 0$ as well, then we see that $(X_{L,T})_G$
  is finitely generated over $\zp$, and we conclude that
  $X_{L,T}$ is finitely generated over $\Lambda(G)$.
\end{proof}

\section{Strongly admissible extensions} \label{general}

In this section, we shall prove our result on the behavior of
inverse limits of minus parts of class groups for
strongly admissible $p$-adic Lie  
extensions of CM-fields
using our extension of Kida's formula (Theorem \ref{kida2}). That
is, we will prove the following theorem, which includes Theorem
\ref{main} (noting Lemma \ref{ven-p-null}).

\begin{theorem} \label{second}
    Let $L/F$ be a strongly admissible $p$-adic Lie extension of
    CM-fields, for an odd prime $p$.
    Set $G = \Gal(L/K)$.
    Let $T$ be a finite set of primes of $F^+$, and let $Q_{L/K}^T$
    be as in \eqref{primeset}.
    Assume that $\mu(X_{K,T}^-) = 0$.
    Then
    $X_{L,T}^-$ is finitely generated over $\Lambda(G)$, and
    we have
    $$
        \rank_{\Lambda({G})} X_{L,T}^- =
        \lambda(X_{K,T}^-)-\delta+|Q_{L/K}^T|.
    $$
\end{theorem}

\begin{remark}
  Let $X_K$ be as in the introduction, and let $Y_K$ be the
  maximal quotient of $X_K$ in which all primes above $p$
  split completely.  Since $\lambda(X_K^-)-\lambda(Y_K^-)$ equals
  the number of primes of $K^+$
  above $p$ which split in $K$ (cf.\ \cite[Proposition 11.4.6]{nsw}),
  Theorem \ref{second} implies that
  $\rk_{\Lambda(G)} X_L^- - \rk_{\Lambda(G)} Y_L^-$
  is the number of primes of $K^+$ above $p$ that split
  completely in $L/K^+$.
  In particular, if no prime above $p$
  splits completely in $L/K^+$, then the kernel of the natural
  surjection from $X_L^-$ to $Y_L^-$ is pseudo-null over $\Lambda(\mc{G})$,
  with $\mc{G} = \Gal(L/F)$.
  This is compatible with \cite[Theorem 4.9]{Ve3}.
  On the other hand, $X_L^-$ is not pseudo-isomorphic to
  $Y_L^-$ if there exists a prime above $p$
  which splits completely in $L/K^+$.
\end{remark}

Let us work in the setting and with the notation of Theorem
\ref{second}. Let $G = \Gal(L/K)$.  We set $d = \dim G$, and for
any prime $v$ of $K$ (or $K^+$), we let $d_v$ denote the dimension
of a decomposition group $G_v$ of $G$ at a prime above $v$. Let
$G_n$ be a sequence of open normal subgroups of $G$ chosen as in
Lemma \ref{well-known}. We then let $g_{n,v}$ denote the index of
	the
image of $G_v$ in the group $G/G_n$.

\begin{lemma} \label{small}
    Let $v$ be a prime of $K$ which does not split completely
    in $L$.  Then $g_{n,v}$ is $O(p^{n(d-1)})$.
\end{lemma}

\begin{proof}
    Let $L_n$ be the subextension in $L/K$ corresponding
    to $G_n$, a finite Galois $p$-extension of $K$.
    Let $C$ be the constant in Lemma \ref{well-known}
    such that $[G:G_n]=Cp^{nd}$
    for all sufficiently large $n$.
    Let $G_v$ and $G_{n,v}$ denote the decomposition groups of
    $G$ and $G_n$, respectively, at a fixed prime of $L$ above
    $v$.  Then, by Lemma \ref{well-known},
    there also exists a rational number $C_v$ such that
    $[G_v:G_{n,v}]=C_v p^{nd_v}$
    for all sufficiently large $n$.
    Since $d_v \ge 1$ by assumption on $v$, we have
    $$
      g_{n,v}=[G/G_n:G_v/G_{n,v}]=(C/C_v)p^{n(d-d_v)}
      =O(p^{n(d-1)}).
    $$
\end{proof}

We are now ready to prove Theorem \ref{second}.

\begin{proof}[Proof of Theorem \ref{second}]
    Again, let $L_n$ be the subextension in $L/K$ corresponding
    to $G_n$.  We remark that $Q_{L_n/K}^T=Q_{L/K}^T$
    if $n$ is sufficiently large.
    By Lemma \ref{small}, $g_{n,v}$ is $O(p^{n(d-1)})$ for
    $v \in Q_{L/K}^T$, so Theorem \ref{kida2} yields
    \begin{equation}\label{estimate}
        \lambda(X_{L_n,T}^-)=
        (\lambda(X_{K,T}^-)-\delta+|Q_{L/K}^T|)[L_n:K]+O(p^{n(d-1)})
    \end{equation}
    for all sufficiently large $n$.

    We let $S_E$ denote the set of primes of an algebraic extension
    $E/K$ consisting of all primes above $p$ and
    the primes which ramify in $L/K$.  We often abbreviate
    $S_E$ simply by $S$.
    Let $L_{n,w}$ denote the completion of $L_n$ at a given prime
    $w$, and let $I_{L_{n,w}}$ denote the inertia group of
    the absolute Galois group $G_{L_{n,w}}$.
    We set
    $$
        \mc{H}_{L_n,T} = \bigoplus_{\substack{w \in S_{L_n}\\w \in T_{L_n}}}
        H^1(G_{L_{n,w}},\qp/\zp) \oplus
        \bigoplus_{\substack{w \in S_{L_n}\\w \notin T_{L_n}}}
        H^1(I_{L_{n,w}},\qp/\zp)
    $$
    and take $\mc{H}_{L,T} = \varinjlim \mc{H}_{L_n,T}$.
    Consider the following commutative diagram.
    $$ \SelectTips{cm}{}
      \xymatrix{
      0 \ar[r] & \Hom({X_{L_n,T}^-},\qp/\zp) \ar[r] \ar[d] &
      H^1(G_{L_n,S},\qp/\zp)^-
      \ar[r] \ar[d]^{\beta_n^-} & \mc{H}_{L_n,T}^-  \ar[d]^{\rho_n^-}\\
      0 \ar[r] & \Hom({X_{L,T}^-},\qp/\zp)^{G_n}  \ar[r] &
      (H^1(G_{L,S},\qp/\zp)^-)^{G_n}
      \ar[r] & (\mc{H}_{L,T}^-)^{G_n}.}
    $$
    Applying the Hochschild-Serre spectral sequence, we
    see that
    $$
        \ker(\beta_n^-)\cong H^1(G_n,\qp/\zp)^-
    $$
    and
    $\coker(\beta_n^-)$ injects into $H^2(G_n,\qp/\zp)^-$.
    By \cite[Lemma 2.5.1]{Ha1},
    the $\zp$-coranks of both of these cohomology groups are bounded by
    a constant $C$, which depends only on $\dim G$, as $n$ increases.
    Applying the snake lemma (with
    $\mc{H}_{L_n,T}^-$ replaced by the image of
    $H^1(G_{L_n,S},\qp/\zp)^-$ in it),
    we have
    \begin{equation}\label{b}
        \lambda(X_{L_n,T}^-)-C \leq \rk_{\Z_p} 
		(X_{L,T}^-)_{G_n}
        \leq \lambda(X_{L_n,T}^-)+\cork_{\Z_p} (\ker(\rho_n^-))+C,
    \end{equation}
    for all sufficiently large $n$.

    We remark that
    \begin{equation*}
        \ker(\rho_n^-)\cong \bigl( \bigoplus_{\substack{w \in S_{L_n} \\ 
        w \notin T_{L_n}}}
        H^1(I_{n,w},\qp/\zp) \oplus \bigoplus_{\substack{w \in S_{L_n} \\ w \in T_{L_n}}}
        H^1(G_{n,w},\qp/\zp) \bigr)^-,
    \end{equation*}
    where $I_{n,w}$ denotes the inertia group of $G_n$ at a prime
    above $w$ (and $G_{n,w}$ is the decomposition group, as
    before).
    The latter equations break up as the direct sum of minus-parts
    of cohomology groups over elements of conjugacy classes of
    primes in $S_{L_n}$ under complex conjugation.  If $w \in S_{L_n}$
    is self-conjugate (or if $w$ splits completely in $L/L_n$),
    then $H^1(I_{n,w},\qp/\zp)^-$ and $H^1(G_{n,w},\qp/\zp)^-$ are trivial.
    If $w$ is complex conjugate to a distinct prime
    $\bar{w}$, then
    $$
      ( H^1(I_{n,w},\qp/\zp) \oplus H^1(I_{n,\bar{w}},\qp/\zp) )^-
      \cong H^1(I_{n,w},\qp/\zp),
    $$
    and similarly with $I_{n,w}$ replaced by $G_{n,w}$.
    We conclude that
    \begin{multline}\label{c}
      \cork_{\zp}(\ker(\rho_n^-)) = \sum_{w \in Q_{L/L_n}^S - Q_{L/L_n}^T}
      \cork_{\zp}H^1(I_{n,w},\qp/\zp) + \\ \sum_{w \in Q_{L/L_n}^T}
      \cork_{\zp}H^1(G_{n,w},\qp/\zp),
    \end{multline}

    The $\zp$-corank of $H^1(G_{n,w},\qp/\zp)$ is less than or equal to
    the
    ${\bf F}_p$-dimension of $H^1(G_{n,w},\Z/p\Z)$, which is
    eventually constant by Lemma \ref{well-known}, equal to the dimension of
    $G_w$.
    Similarly, Lemma \ref{well-known} implies that
    $\cork_{\zp} H^1(I_{n,w},\qp/\zp)$
    is eventually less than or
    equal to the dimension of the inertia group in
    $G_w$, which is at most that of $G_w$.
    Noting that $v\in Q_{L/K}^T$ if and only if $w\in Q_{L/L_n}^T$
    (for any $T$ and)
    for any $w|v$ because $G$ has no elements of finite order,
    equation \eqref{c} implies that
    \begin{equation}\label{d}
      \cork_{\zp}(\ker(\rho_n^-)) \le \sum_{v \in Q_{L/K}^S} d_v g_{n,v}.
    \end{equation}
    Since Lemma \ref{small} implies that
    $g_{n,v}$ is $O(p^{n(d-1)})$ for each
    $v \in Q_{L/K}^S$,
    equations \eqref{b}
    and \eqref{d} yield
    that
    \begin{equation} \label{final}
        \rk_{\zp}((X_{L,T}^-)_{G_n}) = \rk_{\zp} X_{L_n,T}^-
        + O(p^{n(d-1)}) \pm O(1)
    \end{equation}
    As in Lemma \ref{torsion}, the fact that $\mu(X_{K,T}^-) =
    0$ implies that $X_{L,T}^-$ is finitely generated over
    $\Lambda(G)$.
    The result on ranks now follows from
    the fact that $d \ge 2$, equations \eqref{estimate} and
    \eqref{final}, and Lemma \ref{susan}.
\end{proof}

\section{Examples and Remarks} \label{examples}

We conclude this article with several remarks and examples of the
application of Theorem \ref{second}.  We begin by giving three
examples of the application of Theorem \ref{second} to finding $L$
for which $X_L$ is not pseudo-null.

\begin{example} \label{basicex}
  Let $F = \Q(\mu_p)$ and $K = \Q(\mu_{p^{\infty}})$ for an odd prime
  $p$.  Let $\mf{X}_K$ be the Galois group of the maximal pro-$p$ abelian
  unramified
  outside $p$ extension of $K$.  Kummer
  duality induces a pseudo-isomorphism
  $$
    \mf{X}_K^+ \to \Hom_{\zp}(X_K^-,\zp(1))
  $$
  of torsion $\Lambda$-modules with no $\zp$-torsion (see
  \cite[Corollary 11.4.4]{nsw}).
  In particular, we have $\lambda(\mf{X}_K^+) = \lambda(X_K^-)$.
  Thus, there exist (at least) $\lambda(X_K^-)$ distinct
  $\zp$-extensions $L$ of $K$ which are CM, unramified outside $p$, and
  Galois over $\Q$.  For each of these extensions, Theorem
	\ref{second}
  implies that $\rk_{\Lambda(G)} X_L^-
  = \lambda(X_K^-)-1$. In particular, all such
  $X_L^-$ are not pseudo-null if $\lambda(X_K^-) \ge 2$.
\end{example}

\begin{example}
Let $F = \Q(\mu_p)$ and $K=\Q(\mu_{p^\infty})$ for an odd prime
$p$.  We demonstrate how one may construct a $\zp$-extension $L$
of $K$ with $L/F$ Galois such that $X_L^-$ has any sufficiently
large $\Lambda(G)$-rank, for $G = \Gal(L/K)$.

Let $\Pi$ be the set of all rational prime numbers which are
completely decomposed in $F$ but inert in $\Q(\mu_{p^2})/F$. By
the \v{C}ebotarev density theorem, $\Pi$ is a infinite set. Let
$T^+$ denote a finite set of primes of $K^+$ lying above primes in
$\Pi$. Let $S$ be the set of all primes of $K$ above $p$ or a
prime in $T^+$.  From 
	\cite[p.\ 276]{Iw}
we have
\begin{equation}\label{exam}
    0\rightarrow \varinjlim_n X_{K_n}^- \rightarrow H^1(G_{K,S},
    \mu_{p^\infty})^- \rightarrow \bigoplus_{v\in T^+}
    \Q_p/\Z_p \rightarrow 0.
\end{equation}
Furthermore, the summand in the third term which corresponds to
$v$ is canonically the cohomology group of the inertia subgroup at
$v \in T^+$ in $\mf{X}_{K,S}^+$, where $\mf{X}_{K,S}$ is the
Galois group of the maximal abelian pro-$p$ extension of $K$
unramified outside $S$.
Since $v$ splits in $K/K^+$ and lies over an inert prime of $F^+$,
equation \eqref{exam} therefore implies that $\mf{X}_{K,S}^+$
contains a $\Lambda$-submodule isomorphic to
$\Z_p(1)^{\oplus|T^+|}$, with each $\zp(1)$-summand the inertia
group at some $v \in T^+$. Thus, there exists a quotient of
$\mf{X}_{K,S}^+$ by a $\Lambda(\Gamma)$-submodule which defines a
$\Z_p$-extension $L$ that is Galois over $F$, abelian over $K^+$
(hence $L$ is CM), ramified at all primes in $T^+$, and unramified
outside $S$. For this $L$, Theorem \ref{second} yields that
$$
  \rk_{\Lambda(G)}X_L^-=\lambda(X_K^-)+|T^+|-1.
$$
Note that $|T^+|$ can be taken to be arbitrarily large.
\end{example}

\begin{example}
In \cite{Ra}, R.\ Ramakrishna constructs a totally real field $L'$
which is Galois over $\Q$ with Galois group isomorphic to
$PSL_2(\Z_3)$ and which is ramified only at $3$ and $349$. Then,
letting $L = L'\Q(\mu_{3^\infty})$, the Galois group of $L/\Q$ is
a $3$-adic Lie group of dimension four.  It is possible to choose
a number field $F$ contained in $L$ such that $L/F$ is a strongly
admissible $3$-adic Lie extension. Let $K$ denote the cyclotomic
$\Z_3$-extension of $F$.  Applying Theorem \ref{second} to the
extension $L/F$, we see that if $\mu(X_K^-) = 0$, then $X_L^-$ is
not pseudo-null.  However, we do not know that $\mu(X_K^-) = 0$
for this $K$, as $F$ cannot be taken to be a $3$-extension of an
abelian extension of $\Q$.
\end{example}

\begin{remark}
Removing the assumption that $L$ contains the cyclotomic $\zp$-extension
of $K$, Greenberg has recently constructed
nonabelian $p$-adic Lie extensions $L/F$ for which
$\mc{G}=\Gal(L/F)$ is isomorphic to, for instance,
an open subgroup of $PGL_2(\zp)$ and
$X_L$ has nontrivial $\mu$-invariant as a $\Lambda(\mc{G})$-module
(unpublished).
\end{remark}

\begin{remark}
There is an analogous theory for elliptic curves.
Let $E$ be an elliptic curve over a number field $F$, and let $K$ be
the cyclotomic $\zp$-extension of $F$.
In ``classical'' Iwasawa theory for the elliptic curve $E$,
one studies the Pontryagin dual $\mathrm{Sel}_{p^\infty}(E/K)^\vee$
of the Selmer group of $E$ over $K$.
An analogue of Kida's formula for such
Selmer groups is given in \cite{HM} under the assumption that $E$ has good
ordinary reduction at $p$.
One can give a formula for the $\Lambda(G)$-rank of
$\mathrm{Sel}_{p^\infty}(E/L)^\vee$ for a pro-$p$ $p$-adic Lie
extension of $F$ containing $K$ in a similar manner to that of
Theorem \ref{second},
using the same method of proof and a Kida-type formula as above.
In fact, such a formula has already been given in some special
cases: see \cite[Corollary 6.10]{CH} and \cite[Theorem 2.8]{Ho2}
for the case that $L=\Q(E_{p^\infty})$ and \cite[Theorem 3.1]{HV}
for the case in which the dimension of $\Gal(L/F)$ is $2$. In
general, the formula, which is due to the first author, is as
follows. 
We remark that the same formula is also obtained in \cite{bha} 
independently.
\end{remark}

\begin{theorem}
Let $L/F$ be a strongly admissible $p$-adic Lie extension.  Let
$E$ be an elliptic curve defined over $F$ that has good ordinary
reduction at $p$. Let $M_0(L/K)$ be the set of primes $v$ of $K$
not lying above $p$ and ramified in $L/K$, and set
\begin{eqnarray*}
M_1(L/K) &=&\{v\in M_0(L/K): v
     \mr{\ has\ split\ multiplicative\ reduction} \} \\
    M_2(L/K) &=& \{ v\in M_0(L/K) : v \mr{\ has\ good\ reduction\ and\ }
     E(K_v)[p]\ne 0\}.
\end{eqnarray*}
Assume that $\mathrm{Sel}_{p^\infty}(E/K)^\vee$ is finitely
generated over $\Z_p$. Then $\mathrm{Sel}_{p^\infty}(E/L)^\vee$ is
finitely generated over $\Lambda(G)$, and
$$
\rk_{\Lambda(G)} \mathrm{Sel}_{p^\infty}(E/L)^\vee =\rk_{\Z_p}
\mathrm{Sel}_{p^\infty}(E/K)^\vee+ |M_1(L/K)|+ 2|M_2(L/K)|.
$$
\end{theorem}

\appendix
\section{Iwasawa Modules in Procyclic Extensions\\by Romyar T. Sharifi}
\label{special}

In this appendix, we shall derive two types of exact sequences which describe
the behavior of Iwasawa modules in cyclic $p$-extensions $L/K$
such that $L$ is Galois over a number field $F$.  These sequences
and the proof given here are related to a $6$-term exact sequence
of Iwasawa and its method of proof in \cite{iwasawa}, though
derived independently.  By focusing on the case of cyclic
extensions of number fields, we are able to obtain a finer result
than the sequence of Iwasawa, which dealt with general Galois
extensions.  (We also remark that our sequences also bear a relationship with the classical ambiguous class number formula, as in \cite[Lemma 13.4.1]{lang}.)  
We take inverse limits to obtain related sequences for Iwasawa modules in the 
general (pro)cyclic case.

We must first introduce a considerable amount of notation.  For
now, let $K$ be a Galois extension of a number field $F$. In this
section, we allow $p$ to be any prime number. Let $L$ be a cyclic
$p$-extension of $K$ which is Galois over $F$. Set $G =
\Gal(L/K)$, $\mc{G} = \Gal(L/F)$, and $H = \Gal(K/F)$.

Let $T$ be any finite set of primes of $F$ which includes its real
places.  For any algebraic extension $E$ of $F$, let $T_E$ denote
the set of primes of $E$ lying above those in $T$. Let
$\mc{U}_{E,T}$ denote the inverse limit of the $p$-completions of
the $T$-unit groups of the finite subextensions of $F$ in $E$. Let
$X_{E,T}$ denote the maximal unramified abelian pro-$p$ extension
of $E$ in which all primes in $T_E$ split completely.  Let
$\phi_{L/K}^T \colon X_{K,T} \to X_{L,T}^G$ denote the natural
map. Let $\hat{H}^i(M)$ denote the $i$th Tate cohomology group for
the group $G$ and a $\Z[G]$-module $M$.

For each prime $v$ of $K$, let $G_v$ denote the decomposition
group at any prime above $v$ in $G$, and let $I_v$ denote the
inertia group. (If $v$ is a real prime, then 
	$I_v = G_v$
has order $1$ or $2$.) We can consider the map $\Sigma_{L/K}^T$
given by the inverse limit via restriction maps of the sum of
inclusion maps
$$
  \Sigma_{L'/K'}^T \colon \bigoplus_{u \in T_{K'}} G'_u \oplus
  \bigoplus_{u \notin T_{K'}} I'_u \to G',
$$
over number fields $L'$ containing $F$ inside $L$, with $K' = L'
\cap K$, such that $G'_u$ is the decomposition group of
$G'=\Gal(L'/K')$ at $u$ and $I'_u$ is the inertia group.
Similarly, letting $I_{L',T}$ (resp., $I_{K',T}$) denote the
$T$-ideal class group of $L'$ (resp., $K'$) and letting $P_{K',T}$
denote the group of principal $T$-ideals of $K'$, we set
$$
    \mf{I}_{L/K}^T = \lim_{\leftarrow}\,
    ((I_{L',T}^{G'}/I_{K',T}) \otimes_{\Z} \zp)
$$
and
$$    \mf{E}_{L/K}^T = \lim_{\leftarrow}\,
    ((I_{L',T}^{G'}/(P_{K',T} \cdot N_{G'}I_{L',T})) \otimes_{\Z}
    \zp),
$$
where $N_{G'}$ denotes the norm element in $\zp[G']$ and the
inverse limits are taken with respect to norm maps.

Finally, given $r \in \Z$ and an exact sequence of groups
$$
    \Phi \colon\ \ \ldots \to A_i \to A_{i+1} \to \ldots
$$
with a distingushed term $A_0$, let
$$
    \Phi[r] \colon\ \ \ldots \to B_i \to B_{i+1} \to \ldots
$$
denote the exact sequence with distinguished term $B_0$ and $B_i =
A_{i-r}$.

We then have the following theorem.

\begin{theorem} \label{comparison}
    Let $F$ be a number field, $K/F$ a Galois extension, and
    $L/K$ a cyclic $p$-extension with $L/F$ Galois.  Set
    $G = \Gal(L/K)$ and $H = \Gal(K/F)$.  Let $T$ be a finite set
    of primes of $F$.  Then we have canonical exact sequences of
    $\Lambda(H)$-modules:
    \begin{multline*}
    \Gamma_{L/K}^T:\ \
    0 \to (\ker \phi_{L/K}^T) \otimes_{\zp} G \to
    \hat{H}^{-1}(\mc{U}_{L,T}) \to
    \mf{I}_{L/K}^T \otimes_{\zp} G
    \to (\coker \phi_{L/K}^T) \otimes_{\zp} G
    \to \\ \hat{H}^0(\mc{U}_{L,T}) \to
    \ker \Sigma_{L/K}^T \to (X_{L,T})_G \to
    X_{K,T} \to \coker \Sigma_{L/K}^T \to 0
    \end{multline*}
    and
    \begin{multline*}
    \Psi_{L/K}^T:\ \
    \ldots \to \hat{H}^{-1}(\mc{U}_{L,T}) \to \mf{E}_{L/K}^T \otimes_{\zp} G
    \to \hat{H}^0(X_{L,T}) \otimes_{\zp} G \to \\ \hat{H}^0(\mc{U}_{L,T}) \to
    \ker \Sigma_{L/K}^T \to \hat{H}^{-1}(X_{L,T}) \to \ldots
    \end{multline*}
    with $\Psi_{L/K}^T[6] = \Psi_{L/K}^T \otimes_{\zp} G$.
\end{theorem}

\begin{proof}
    We will leave out subscript and superscript $T$'s throughout
    this proof, for compactness.
    For $E$ a finite extension of $F$, we let $\mc{O}_{E}$ denote
    the ring of $T$-integers of $E$, let $I_{E}$ denote the group of
    fractional ideals of $\mc{O}_{E}$, let $P_{E}$ denote the
    subgroup of principal fractional ideals, and let $\Cl_{E}$
    denote the class group of $\mc{O}_{E}$.

    We begin by proving the theorem in the case that $K$ and $L$ are both
    number fields.  
    Consider the commutative diagram
    of exact sequences:
    \begin{equation} \label{firstseq} \SelectTips{cm}{}
        \xymatrix@R=10pt@C=5pt{
        &&&& 0 \ar[rd] &&0 \\
        &&&&& P_E \ar[ur] \ar[rd] \\
        0 \ar[rr] && \UT{E} \ar[rr] && E^{\times} \ar[rr] \ar[ur] &&
        I_E \ar[rr] && \Cl_E \ar[rr] && 0.
        }
    \end{equation}

    We obtain from \eqref{firstseq} in the case $E = L$ two long exact
    sequences in Tate cohomology
    \begin{equation} \label{twoseqs} \small \SelectTips{cm}{}
        \xymatrix@C=6pt{
        & 0 \ar[r] & \hat{H}^{2i-1}(P_{L}) \ar[r] \ar@{=}[d] &
        \hat{H}^{2i}(\UT{L}) \ar[r] & \hat{H}^{2i}(L^{\times}) \ar[r] &
        \hat{H}^{2i}(P_{L}) \ar[r] \ar@{=}[d] & \hat{H}^{2i+1}(\UT{L}) \ar[r]
        & 0 \\
        \ldots \ar[r] & \hat{H}^{2i-2}(\Cl_L) \ar[r] &
        \hat{H}^{2i-1}(P_{L}) \ar[r] &
        0 \ar[r] & \hat{H}^{2i-1}(\Cl_L) \ar[r] &
        \hat{H}^{2i}(P_{L}) \ar[r] & \hat{H}^{2i}(I_{L}) \ar[r] &
        \ldots,}
    \end{equation}
    where we have used that $\hat{H}^{2i-1}(L^{\times}) =
    \hat{H}^{2i-1}(I_L) = 0$.
    Chasing the diagram \eqref{twoseqs}, we obtain an exact
    sequence:
    \begin{multline} \label{longseq}
        \ldots \to \hat{H}^{2i-2}(\Cl_{L}) \to
        \hat{H}^{2i}(\UT{L}) \to
        \ker(\hat{H}^{2i}(L^{\times}) \to \hat{H}^{2i}(I_{L})) \to \\
        \hat{H}^{2i-1}(\Cl_{L}) \to \hat{H}^{2i+1}(\UT{L})
         \to \coker(\hat{H}^{2i}(L^{\times}) \to \hat{H}^{2i}(I_{L})) \to \ldots.
    \end{multline}

  Using \eqref{firstseq} again, we have a commutative diagram
  \begin{equation*} \label{seconddiag} \SelectTips{cm}{}
    \xymatrix@!C=24pt{
    0 \ar[r] & \UT{K} \ar[r] \ar@{=}[d] & K^{\times} \ar[r] \ar@{=}[d] &
    P_K \ar[r] \ar[d] & 0 \\
    0 \ar[r] & \UT{K} \ar[r] & K^{\times} \ar[r] & P_L^G \ar[r] &
    \hat{H}^1(\UT{L}) \ar[r] & 0,
    }
  \end{equation*}
  and it provides an isomorphism $\hat{H}^1(\UT{L}) \cong P_L^G/P_K$.
  Noting this and applying the snake lemma to the commutative diagram
  \begin{equation*} \label{firstdiag} \SelectTips{cm}{}
    \xymatrix@!C=24pt{
    0 \ar[r] & P_K \ar[r] \ar[d] & I_K \ar[r] \ar[d] &
    \Cl_K \ar[r] \ar[d]^{\phi_{L/K}} & 0 \\
    0 \ar[r] & P_L^G \ar[r] & I_L^G \ar[r] & \Cl_L^G \ar[r] &
    \hat{H}^1(P_L) \ar[r] & 0,
    }
  \end{equation*}
  we obtain an exact sequence
  \begin{equation} \label{secondseq}
    0 \to \ker(\Cl_K \xrightarrow{\phi_{L/K}} \Cl_L) \to \hat{H}^1(\UT{L}) \to
    I_L^G/I_K \to \Cl_L^G/\phi_{L/K}(\Cl_K) \to \hat{H}^1(P_L) \to 0.
  \end{equation}

  One next checks easily that the map
  $\hat{H}^{-1}(\Cl_L) \to \hat{H}^1(\mc{O}_L^\times)$ in \eqref{longseq}
  has image contained in $\ker(\phi_{L/K})$ via the map in
  \eqref{secondseq} and that the resulting map $\hat{H}^{-1}(\Cl_L) \to \Cl_K$
  is induced by the norm map $(\Cl_L)_G \to \Cl_K$.  Furthermore, since
  the kernel of the norm map is contained in $\hat{H}^{-1}(\Cl_L)$, we have
  an exact sequence
  \begin{equation} \label{midway}
    \ldots \to \hat{H}^{-2}(\Cl_L) \to
    \hat{H}^0(\UT{L}) \to \ker(\hat{H}^0(L^{\times}) \to
    \hat{H}^0(I_L)) \to (\Cl_L)_G \to \Cl_K.
  \end{equation}
  Next, we attach \eqref{secondseq} to the left of
  \eqref{midway} via the map $\hat{H}^{-1}(P_L) \to
  \hat{H}^0(\UT{L})$ in \eqref{twoseqs}, obtaining
  \begin{multline} \label{almost} \SelectTips{cm}{}
    0 \to \ker \phi_{L/K} \otimes G \to
    \hat{H}^{-1}(\UT{L}) \to I_L^G/\phi(I_K) \otimes G \to
    \coker \phi_{L/K} \otimes G \to \\
    \hat{H}^0(\UT{L}) \to  \ker(\hat{H}^0(L^{\times}) \to
    \hat{H}^0(I_L)) \to  (\Cl_L)_G \to \Cl_K.
  \end{multline}

  Now, we must study the map
  $\hat{H}^0(L^{\times}) \to \hat{H}^0(I_L)$.  By class field theory,
  we have an exact sequence
  $$
    0 \to \hat{H}^2(L^{\times}) \to
    \hat{H}^2\bigl(\bigoplus_w L_w^{\times}\bigr) \to
    \frac{1}{[L:K]}\Z/\Z,
  $$
  in which $L_w$ denotes the completion of $L$ at a prime $w$.
  This yields a sequence that fits into a commutative diagram
  \begin{equation} \label{anotherdiag} \SelectTips{cm}{}
    \xymatrix{
    0 \ar[r] & \hat{H}^0(L^{\times}) \ar[d] \ar[r] &
    \hat{H}^0(\bigoplus_w L_w^{\times})
    \ar[r] \ar@{->>}[d] & G \\
    & \hat{H}^0(I_L) & \hat{H}^0(\bigoplus_{w \notin T_L} L_w^{\times})
    \ar[l] &
    \hat{H}^0(\bigoplus_{w \notin T_L} U_w^{\times}) \ar[l] & 0 \ar[l],
    }
  \end{equation}
  where $U_w$ denotes the unit group of $L_w$.
  The local reciprocity maps provide canonical isomorphisms
  \begin{eqnarray*}
    \hat{H}^0\bigl(\bigoplus_{w \mid v} L_w^{\times}\bigr) \cong
    G_v & \mr{and} &
    \hat{H}^0\bigl(\bigoplus_{w \mid v} U_w^{\times}\bigr) \cong I_v
  \end{eqnarray*}
  for any prime $v$ of $K$.
  We also have a non-canonical isomorphism
  \begin{equation} \label{idealcohom}
    \hat{H}^0(I_L) \cong \bigoplus_{v \notin T_K} G_v.
  \end{equation}
  Making the resulting replacements in
  \eqref{anotherdiag}, we have a commutative diagram
  $$ \SelectTips{cm}{}
    \xymatrix{
    0 \ar[r] & \hat{H}^0(L^{\times}) \ar[d] \ar[r] &
    \bigoplus_v G_v
    \ar[r] \ar@{->>}[d] & G \\
    &  \bigoplus_{v \notin T_K} G_v & \bigoplus_{v \notin T_K} G_v
    \ar_{\oplus_v |I_v|}[l] & \bigoplus_{v \notin T_K} I_v \ar[l] & 0, \ar[l]
    }
  $$
  which implies that we have a canonical isomorphism
  \begin{equation} \label{kerident}
    \ker(\hat{H}^0(L^{\times}) \to \hat{H}^0(I_L)) \cong \ker \Sigma_{L/K}.
  \end{equation}
  Furthermore, we remark that
  \begin{equation} \label{cokerident}
    \coker( \hat{H}^0(L^{\times}) \to \hat{H}^0(I_L) ) \cong
    I_L^G/(P_K \cdot N_G I_L).
  \end{equation}

  Plugging \eqref{kerident}
  into \eqref{almost} and noting that
  $$
    \coker(\Cl_L \to \Cl_K) \cong \coker \Sigma_{L/K},
  $$
  we obtain an exact sequence with desired $p$-part $\Gamma_{L/K}^T$.
  Similarly, plugging
  \eqref{kerident} and \eqref{cokerident}
  into \eqref{longseq}, we obtain the sequence $\Psi_{L/K}^T$.

  Note that there are natural maps of exact sequences, $\Gamma_{L'/K'}^T \to
  \Gamma_{L/K}^T$, for finite extensions $L'/L$ and $K'/K$ such that
  $L'/K'$ is cyclic,
  which are given by norm maps from $L'$ to $L$
  on the first, second, third, fourth, and seventh terms,
  norm maps from $K'$ to $K$ on the fifth and eighth terms,
  and the maps induced by restriction of Galois groups on the
  remaining two terms.
  The sequence $\Gamma_{L/K}^T$ 
	in the  
  case in which $L$ is not a number field
  now follows by taking the inverse limit of
  the sequences $\Gamma_{L'/L' \cap K}^T$, with $L'$ a number field
  contained in $L$.
  Since $T$ is assumed to be finite,
  all terms at the finite level are finite, and therefore, the sequences
  remain exact in the inverse limit.
  Similarly, we
  have maps $\Psi_{L'/K'}^T \to \Psi_{L/K}^T$, and the inverse
  limit yields the desired sequence in the general case.
\end{proof}

Taking inverse limits, one easily obtains the following corollary
for $\zp$-extensions $L/K$ with Galois group $G$.  Here, we let
$N_{L/K}^T$ denote the obvious map $(\mc{U}_{L,T})_G \to
\mc{U}_{K,T}$ induced by the inverse limit of norm maps.

\begin{corollary} \label{zpext}
    Let $F$ be a number field, $K/F$ a Galois extension, and
    $L/K$ a $\zp$-extension with $L/F$ Galois.  Set
    $G = \Gal(L/K)$ and $H = \Gal(K/F)$.  Let $T$ be a finite set
    of primes of $F$.  Then we have canonical exact sequences of
    $\Lambda(H)$-modules:
    \begin{multline*}
    \Gamma_{L/K}^T:\ \
    0 \to \ker N_{L/K}^T \to
    \mf{I}_{L/K}^T \otimes_{\zp} G
    \to X_{L,T}^G \otimes_{\zp} G
    \to \coker N_{L/K}^T \to \\
    \ker \Sigma_{L/K}^T \to (X_{L,T})_G \to
    X_{K,T} \to \coker \Sigma_{L/K}^T \to 0
    \end{multline*}
    and
    \begin{multline*}
    \Psi_{L/K}^T:\ \ \ldots
    \to \ker N_{L/K}^T \to \mf{E}_{L/K}^T \otimes_{\zp} G
    \to X_{L,T}^G \otimes_{\zp} G \to \\
    \coker N_{L/K}^T \to \ker \Sigma_{L/K}^T \to (X_{L,T})_G \to
    \ldots
    \end{multline*}
    with $\Psi_{L/K}^T[6] \cong \Psi_{L/K}^T \otimes_{\zp} G$.
\end{corollary}

Corollary \ref{zpext} can be used to give a simple proof of
Theorem \ref{second} for $\zp$-extensions which includes the case
of $p = 2$, though we do not include it here.  There are numerous
remarks to be made.

\begin{remarks}
\ \\ \vspace{-4ex}
  \begin{enumerate}
  \item[1.]  Let $R_{L/K}$ denote the
  set of finite primes $v$ of $K$ with $v \notin T_K$ and such that
  the completion $L_w$ for $w \mid v$ does not contain the unramified
  $\zp$-extension of $F_v$.
  We have a noncanonical isomorphism of $\Lambda(H)$-modules
  $$
    \mf{I}_{L/K}^T \cong
    \prod_{v \in R_{L/K}} I_v
  $$
  and a canonical exact sequence
  $$
    0 \to \coker \Sigma_{L/K}^T \to \mf{E}_{L/K}^T \to
    \mf{I}_{L/K}^T \to 0.
  $$
  Note that, if $G \cong \zp$, then
  the set of $v \in R_{L/K}$ with
  $I_v \neq 0$ consists only of primes over $p$.
  \item[2.] It is not necessary to assume that $G$ is pro-$p$, rather
  just procyclic, if we replace Galois groups by their
  $p$-completions in the definitions of terms of $\Gamma_{L/K}^T$,
  aside from $G$-invariants and coinvariants.
  \item[3.] In addition to the functoriality inducing maps
  $\Gamma_{L'/K'}^T \to \Gamma_{L/K}^T$ and $\Psi_{L'/K'}^T \to
  \Psi_{L/K}^T$, there are natural maps
  $\Gamma_{L/K}^T \to \Gamma_{L/K}^{T'}$ and $\Psi_{L/K}^T \to
  \Psi_{L/K}^{T'}$
  for $T \subseteq T'$ (induced by the natural quotient maps on
  class groups and ideal groups, inclusion maps on unit groups, and
  equality on $G$ together with the natural maps between its subgroups).
  \item[4.] If we only wish to have an exact sequence of
  $\zp$-modules, we need not assume that $L$ and $K$ are Galois
  over a number field $F$.  To do this, choose a set $T_K$ of
  primes of $K$ containing the real places,
  let $T_E$ be the set of primes lying below
  those in $T_K$ for any $E \subset K$, and assume that $T = T_{\Q}$ is
  finite.  The sequences $\Gamma_{L/K}^T$ and $\Psi_{L/K}^T$
  of $\zp$-modules defined as before are still exact.
  \item[5.] It is not necessary to assume that $T$ is finite if $L$ is
  a number field, since we do not have to pass to an inverse
  limit.
  \item[6.] If $T_K$ contains the set of primes which ramify in $L/K$, then
  we have $I_v = 0$ for all $v \notin T_K$, so 
	when $G \cong {\bf Z}_p$  
  we obtain a $6$-term exact sequence
  \begin{multline*}
    \qquad \ \ 0 \to 
	X_{L,T}^G
    \otimes_{\zp} G
    \to \coker N_{L/K}^T \to \\
    \ker \Sigma_{L/K}^T \to (X_{L,T})_G \to
    X_{K,T} \to \coker \Sigma_{L/K}^T \to 0
  \end{multline*}
  from $\Gamma_{L/K}^T$, and $\Psi_{L/K}^T$ becomes
  \begin{multline*}
    \qquad \ \ \ldots \to \ker N_{L/K}^T \to \coker \Sigma_{L/K}^T 
	\otimes_{{\bf Z}_p} G    
    \to X_{L,T}^G \otimes_{\zp} G \to\\ \coker N_{L/K}^T \to
    \ker \Sigma_{L/K}^T \to (X_{L,T})_G \to \ldots.
  \end{multline*}
  \item[7.] If $K$ contains the cyclotomic $\zp$-extension of $F$, then only
  (the $p$-parts of) those $G_v/I_v$ with $v$ lying
  over $p$ (or real places when $p = 2$) can be nontrivial.
  \end{enumerate}
\end{remarks}

We now mention a couple of other approaches to the proof of
Theorem \ref{comparison} for the sequences $\Gamma_{L/K}^T$, as we
believe the methods are quite interesting in their own right and
may apply in other contexts. Perhaps surprisingly, the method we
have given above is not only the most down-to-earth approach but
also seemingly the most easily applied to treat the general case.
In describing the alternate approaches, we focus on the case that
$T$ contains the primes above $p$ and all primes which ramify in
$L/K$, for which Galois cohomology with restricted ramification is
most easily applied.

The first approach involves again working first in the case that
$L$ is a number field and writing out a seven-term exact sequence
using the Hochschild-Serre spectral sequence
$$
    H^p(G,H^q(G_{L,T},\mc{O}_{\Omega,T}^{\times})) \Rightarrow
    H^{p+q}(G_{K,T},\mc{O}_{\Omega,T}^{\times}),
$$
where $G_{E,T}$ denotes the maximal unramified outside $T$
extension of $E$ (for $E = K,L$) and $\mc{O}_{\Omega,T}$ is the
ring of $T$-integers of $\Omega$, the maximal unramified outside
$T$ extension of $K$.  Note that one has nice descriptions of (the
$p$-completions of) the groups
$H^i(G_{K,T},\mc{O}_{\Omega,T}^{\times})$ for $i = 0,1,2$ in terms
of units, class groups, and the Brauer group, respectively. (At
one point, one must describe explicitly
the map $E_2^{1,1} \to E_2^{3,0}$,
and this requires comparing with the map $\hat{H}^{-1}(\Cl_{L,T})
\to \hat{H}^1(\mc{O}_{L,T}^{\times})$ in \eqref{longseq}.) One
then passes to the inverse limit.

Another approach involves using the Poitou-Tate sequences for $E$
equal to $K$ and $L$ of \cite[Theorem 5.4]{jannsen}:
\begin{equation} \label{pt}
    0 \to H^2(G_{E,T},\qp/\zp)^{\vee}
    \to \mc{U}_{E,T} \to \mc{A}_{E,T}
    \to \mf{X}_{E,T} \to X_{E,T} \to 0
\end{equation}
where $\mf{X}_{E,T}$ denotes the Galois group of the maximal
	abelian
pro-$p$ unramified outside $T$ extension of $E$ and
$$
    \mc{A}_{E,T} =
    \lim_{\substack{\leftarrow\\E \subset K}}
    \bigoplus_{w \in T_E} H^1(G_{E_w},\qp/\zp)^{\vee},
$$
where $G_{E_w}$ is the absolute Galois group of the completion
$E_w$ of $E$ at $w$ and ${}^{\vee}$ is used to denote the
Pontryagin dual. Let us focus on the case $G \cong \zp$, which one
can treat directly. In this case, one breaks up the $5$-term exact
sequences \eqref{pt} into three pairs of $3$-term exact sequences
and then considers maps on $G$-coinvariants between them  (the
general idea here being taken from \cite{cs}).  One obtains three
long exact sequences via the snake lemma and derives the desired
$6$-term sequence from these, using repeatedly the fact that $G$
has $p$-cohomological dimension $1$. (In fact, in the case that
$H^2(G_{K,T},\qp/\zp)$ is nontrivial, we only carried it out up to
a certain difficult check of commutativity.) One can also derive
the desired exact sequence from the spectral sequence for a
certain three-by-five complex with exact rows consisting of two
copies of \eqref{pt} for $L$ and one for $K$, which degenerates at
$E_4$, and we thank Marc Nieper-Wisskirchen for suggesting the
idea.

\renewcommand{\baselinestretch}{1}

\vspace{2ex} \footnotesize \noindent
Yoshitaka Hachimori\\
CICMA\\
Department of Mathematics and Statistics\\ 
Concordia University\\
1455 de Maisonneuve Blvd. West\\
Montr\'{e}al, Qu\'{e}bec H3G 1M8, Canada\\
e-mail address: {\tt yhachi@mathstat.concordia.ca}\\ \\
Romyar Sharifi \\
Department of Mathematics and Statistics\\
McMaster University\\
1280 Main Street West\\
Hamilton, Ontario L8S 4K1, Canada\\
e-mail address: {\tt sharifi@math.mcmaster.ca}
\end{document}